\documentstyle[11pt]{article}

\newtheorem{thm}{Theorem}[section]
\newtheorem{lem}[thm]{Lemma}

\newtheorem{pro}[thm]{Proposition}

\newtheorem{defi}[thm]{Definition}

\newcommand{\gm }{\Gamma }
\newcommand{\lon }{\longrightarrow }
\newcommand{\be }{\begin{eqnarray*}}
\newcommand{\ee }{\end{eqnarray*}}

\input amssym.def
\input amssym

\newcommand{\frakg}{{\frak g}}
\newcommand{\frakh}{{\frak h}}

\newcommand{\hstar}{*_{\hbar}}

\newcommand{\cald}{{\cal D}}

\newcommand{\smalcirc}{\mbox{\tiny{$\circ $}}}

\def\description label#1{\hfil\bf[#1]\hfil}
\newcommand{\al}{\alpha}
\newcommand{\lrw}{\longrightarrow}
\newcommand{\ot}{\mbox{$\otimes$}}
\newcommand{\Map}{\longmapsto}

\newcommand{\otr}{\ot_{R}}
\newcommand{\oth}{\ot_{\hbar}}
\newcommand{\otrh}{\ot_{R_{\hbar}}}
\newcommand{\dh}{\cald [[\hbar ]]}
\newcommand{\dhh}{\cald_{\hbar}}
\newcommand{\limh}{lim_{\hbar \mapsto 0}}

\def\sdp{\mathbin{\hbox{$\mapstochar\kern-.3333em\times$}}}
\def\pds{\mathbin{\hbox{$\times\kern-.55em\mapstochar\,$}}}

\newcommand{\wed}{\mathbin{\lower1.5pt\hbox{$\scriptstyle{\wedge}$}}}

\let\Tilde=\widetilde

\def\chigh{{\raise1.5pt\hbox{$\chi$}}}
\let\phi=\varphi
\def\til0{\Tilde{0}}

\def\dminus{\raise2pt\hbox{\vrule height1pt width 2ex}\hskip3pt}

\def\pback#1{\mathbin{{{\lower1.2ex\hbox{$\times$}}\atop #1}}}

\def\vlra{\hbox{$\,-\!\!\!-\!\!\!-\!\!\!-\!\!\!-\!\!\!
-\!\!\!-\!\!\!-\!\!\!-\!\!\!-\!\!\!\longrightarrow\,$}}

\def\gpd{\,\lower1pt\hbox{$\longrightarrow$}\hskip-.24in\raise2pt
             \hbox{$\longrightarrow$}\,}

\def\lgpd{\,\lower1pt\hbox{$\vlra$}\hskip-1.02in\raise2pt\hbox{$\vlra$}\,}

\def\llgpd{\,\lower1pt\hbox{$\vvlra$}\hskip-1.3in\raise2pt\hbox{$\vvlra$}\,}

\begin{document}

\title{{\bf Quantum groupoids and deformation quantization }}

\author{ PING XU \thanks{ Research partially supported by NSF
        grants  DMS95-04913 and DMS97-04391.}\\
Research Institute for Math. Science\\
Kyoto University, Kyoto, JAPAN\\
 AND\\
 Department of Mathematics\\The  Pennsylvania State University \\
University Park, PA 16802, USA\\
        {\sf email: ping@math.psu.edu }}

\date{}

\maketitle
\begin{abstract}
The purpose of this Note is to unify quantum groups and star-products 
 under a general umbrella: quantum groupoids.
It is shown that a quantum groupoid
  naturally  gives rise to a Lie bialgebroid 
as a classical limit.
The converse question, i.e., the quantization problem, 
 is posed.  In particular, any regular triangular Lie bialgebroid
is shown quantizable.   For the Lie bialgebroid of a Poisson manifold,
its   quantization is equivalent to   a star-product.
\end{abstract}

{\bf Groupo\"{\i}des quantiques et quantification par d\'eformation}

{\bf R\'esum\'e}
Cette note a pour but d'unifier groupes quantiques et star-produits sous
une m\^eme enseigne: les {\it groupo\"{\i}des quantiques}. Nous montrons
que tout groupo\"{\i}de quantique admet de mani\`ere naturelle une
big\'ebro\"{\i}de de Lie comme limite classique. Le probl\`eme
r\'eciproque de quantification est pos\'e, et nous le r\'esolvons dans
le cas des big\'ebro\"{\i}des de Lie triangulaires r\'eguli\`eres.
Enfin, quantifier la big\'ebro\"{\i}de associ\'ee \`a une vari\'et\'e de
Poisson revient \`a y construire un star-produit.

{\bf Version fran\c caise abr\'eg\'ee}
Les tenseurs de Poisson, sous plusieurs aspects, 
ressemblent aux $r$-matrices triangulaires classiques de la th\'eorie 
des groupes quantiques. La notion de big\'ebro\"{\i}de de Lie
introduite dans \cite{MackenzieXu1} permet d'unifier les structures de Poisson 
et les big\`ebres de Lie. Le th\'eor\`eme d'int\'egration pour les
big\'ebro\"{\i}des de Lie  contient \`a la fois le th\'eor\`eme 
d'int\'egration de Drinfeld pour les big\`ebres de Lie et le
th\'eor\`eme de Karasev-Weinstein sur l'existence locale
de groupo\"{\i}des symplectiques sur une vari\'et\'e de Poisson.
Les groupes quantiques apparaissent comme quantification
de big\`ebres de Lie, alors que la quantification des strutures
de Poisson se r\'ealise \`a l'aide des star-produits.
Il y a donc tout lieu d'esp\'erer que ces deux objets quantiques
sont intimement li\'es.
Le but de cette note est d'\'etablir un rapport  entre
star-produits et groupes quantiques dans le cadre g\'en\'eral
des groupo\"{\i}des quantiques (ou des QUE-alg\'ebro\"{\i}des), 
c'est-\`a-dire des d\'eformations d'alg\'ebro\"{\i}des de Hopf de 
l'alg\`ebre enveloppante universelle d'une
alg\'ebro\"{\i}de de Lie.

Soit $A$ une alg\'ebro\"{\i}de de Lie, son alg\`ebre enveloppante 
universelle $UA$ poss\`ede une structure d'alg\'ebro\"{\i}de de Hopf
cocommutative. En particulier, lorsque $A$ est l'alg\'ebro\"{\i}de de 
Lie d'un fibr\'e tangent TP, cette structure d'alg\'ebro\"{\i}de de
Hopf est celle de l'alg\`ebre ${\cal D}(P)$ des op\'erateurs
diff\'erentiels
sur $P$.

Les principaux r\'esultats pr\'esent\'es dans cette note sont:\\\\
{\bf Theorem A} \ \ {\em A tout  star-produit sur une vari\'et\`e  $P$ 
correspond un
groupo\"{\i}de quantique $\cald_{\hbar} (P) $,
d\'eformation de l'alg\'ebro\"{\i}de de Hopf de $\cald (P)$. }\\\\
{\bf Theorem B} \ \ {\em Un groupo\"{\i}de quantique  $U_{\hbar} A$ admet de mani\`ere naturelle
une big\'ebro\"{\i}de de Lie $(A, A^* )$ comme limite classique.
La structure de Poisson induite par cette big\'ebro\"{\i}de de Lie
coincide avec celle qui est associe\'e  \`a  la $*$-alg\`ebre 
$C^{\infty}(P)[[\hbar ]]$.}

R\'eciproquement, une quantification d'une big\'ebro\"{\i}de de 
Lie $(A, A^{*})$ est un groupo\"{\i}de quantique $U_{\hbar}A$ 
dont la limite classique est $(A, A^{*})$. \\\\
{\bf Theorem C} \ \ {\em Toute big\'ebro\"{\i}de de Lie triangulaire  r\'eguli\`ere admet
une quantification.}

Nous croyons que, dans le cadre g\'en\'eral d\'ecrit ci-dessus, 
les techniques
utilis\'ees dans la th\'eorie des groupes quantiques  
peuvent donner une meilleure compr\'ehension des
star-produits sur une vari\'et\'e de Poisson.

\section{Introduction}

Poisson  tensors in many aspects  resemble classical triangular r-matrices
in quantum group theory. A  notion unifying both Poisson structures
and Lie bialgebras was introduced in \cite{MackenzieXu1}, called
   Lie bialgebroids.
Integration theorem for Lie bialgebroids encomposes both Drinfeld
theorem of integration of Lie bialgebras on the one hand,
and Karasev-Weinstein theorem of existence  of  local symplectic groupoids
for Poisson manifolds on the other hand. Quantization of
Lie bialgebras leads  to quantum groups, while quantization of Poisson
manifolds is the so called   star-products.
It is therefore natural to expect that there should exist  some intrinsic
connection between these two quantum objects.
The  purpose  of the Note is to connect  these two concepts
under the  general framework of  quantum groupoids (or QUE algebroids), i.e.,
Hopf algebroid deformation of  the universal enveloping algebras of
Lie algebroids.


Given  a Lie algebroid $A$,
its universal enveloping algebra $UA$  carries   a  natural cocommutative 
Hopf algebroid structure.  In particular, when $A$ is the
tangent bundle Lie algebroid $TP$, this 
is the Hopf algebroid structure on   $\cald (P)$, the algebra  of differential
operators on $P$.


We will see that  a star-product on $P$  corresponds to a quantum
groupoid $\cald_{\hbar} (P) $,
 which is a Hopf algebroid  deformation of  $\cald (P)$.
In  general,  we show that
a  quantum groupoid  $U_{\hbar }A$ naturally induces a Lie bialgebroid
$(A, A^* )$ as a classical limit.
Then, we pose  the general question  of    quantization  of
Lie bialgebroids, which, as special
cases,  encomposes quantization of Lie bialgebras
and deformation quantization of Poisson manifolds.
In particular, we   prove that any regular
 triangular Lie bialgebroid is quantizable.

Our main motivation is that  this general framework may provide
some new insights 
 in understanding star  products of Poisson manifolds via
the  methods in quantum group theory.

It is worth noting  that the notion of Hopf algebroids was  introduced
by Lu \cite{lu} essentially by translating the axioms of Poisson groupoids
to their quantum counterparts,  while the case that base algebras
are commutative
 was already
studied  by Maltsiniotis \cite{M}  in 1992.

\section{Hopf algebroids}
 
\begin{defi}
\label{dfn_algebroid}

A  Hopf algebroid consists of the following data:

1) a  {\bf total algebra} $H$ with product $m$, a  base algebra $R$, 
  a  {\bf source map}: an algebra homomorphism $\al: ~ R \lrw H$, 
and  a {\bf target map}: an algebra anti-homomorphism
$\beta: ~ R \lrw H $
such that the images of $\al$ and $\beta$ commute in $H$, i.e., $\forall
a, b \in R$,
$\al(a) \beta (b) ~ = ~ \beta (b) \al (a)$.
There is then a natural $(R ,  R)$-bimodule structure on $H$ given by
 $a \cdot  h = \al (a)h$ and 
$ h \cdot  a = \beta (a) h$.
Thus, we can form the $(R, R)$-bimodule
product $H \otr  H$.  It is easy to see that $H \otr  H$  again admits
an  $(R, R)$-bimodule structure.  This will allow us to form the 
 triple product $H \otr H \otr  H$ and so on.

2) a  {\bf co-product}: an $(R, R)$-bimodule map 
$\Delta: ~ H \lrw H \otr  H$
with $\Delta(1) = 1 \ot 1$ satisfying the co-associativity:
\begin{equation}
\label{eq:coassociative}
(\Delta \otr  id_{\scriptscriptstyle H} ) \Delta 
~ = ~ (id_{\scriptscriptstyle H} \otr  \Delta)
 \Delta: ~ H \longrightarrow H \otr H \otr H;
\end{equation}

3) the product and the co-product  are  
{\bf compatible} in the following sense:
\begin{equation}
\label{eq:compatible1}
\Delta (h)(\beta (a)\ot 1-1\ot \alpha (a))=0,  \  \  \mbox{ in }
H\otr H, \ \forall a\in R \mbox{ and }
h\in H, \mbox{ and }
\end{equation}

\begin{equation}
\label{eq:compatible2}
\Delta (h_{1}h_{2})=\Delta (h_{1})\Delta (h_{2}), \ \ \forall h_{1}, h_{2}\in H,
\ \ \ \ \mbox{(see the  remark below)};
\end{equation}

4) a  co-unit map: an (R, R)-bimodule map
$\epsilon: ~~ H \lrw R$
satisfying 
$\epsilon(1_{\scriptscriptstyle H}) ~ = ~ 1_{\scriptscriptstyle R} $
(it follows then that $\epsilon \beta   =  \epsilon \al  =  id_{\scriptscriptstyle R}$) and
\begin{equation}
\label{eq_co-unit}
 (\epsilon \otr id_{\scriptscriptstyle H}) \Delta ~ = ~ 
(id_{\scriptscriptstyle H} \otr \epsilon) \Delta ~ = 
~ id_{\scriptscriptstyle H}: ~ H \lrw H.
\end{equation}
Here  we have used the identification:  $R\otr H\cong H\otr R\cong H$
(note that both  maps on the left hand sides of Equation  (\ref{eq_co-unit})
are well-defined).

We will denote this   Hopf algebroid by $(H, R, \alpha ,
 \beta , m , \Delta , \epsilon)$.

\end{defi}
{\bf Remark} (1) It is clear that any left $H$-module is automatically
an $(R, R)$-bimodule.
Now given  any  left $H$-modules $M_{1}$ and $M_{2}$,  for $m_{1}\in M_{1}$,
$m_{2}\in M_{2}$ and $h\in H$,
define,
\begin{equation}
h\cdot (m_{1}\otr m_{2})=\Delta (h)(m_{1}\ot m_{2}).
\end{equation}
The RHS is a well-defined element in $M_{1}\otr M_{2}$ due to Equation (\ref{eq:compatible1}).
In particular, when  taking  $M_{1}=M_{2}=H$, we see that the RHS of
Equation (\ref{eq:compatible2}) makes sense.
In fact,   Equation (\ref{eq:compatible2}) implies
 that $M_{1}\otr M_{2}$ is again
a left $H$-module.

(2)  The compatibility condition (Equations 
(\ref{eq:compatible1})  and (\ref{eq:compatible2}) )
is   equivalent to the following one  in Lu's original  definition \cite{lu}:
 the kernel of the  map
\begin{equation}
\label{eq_phi}
\Psi:~
H \ot H \ot H \lrw H \otr H: ~~  \sum h_1 \ot h_2 \ot h_3
 \Map \sum (\Delta h_1) (h_2 \ot h_3)
\end{equation}
is a  left ideal of $H \ot H^{op} \ot H^{op}$, where $H^{op}$
denotes $H$ with the opposite product. 


(3) 
In our definition above,  we choose  not to require  the existence of
 antipodes  because many 
 interesting examples, as shown   below, 
often  do not admit antipodes.

The following proposition follows immediately from definition.

\begin{pro}
Let $(H, R, \alpha , \beta , m , \Delta , \epsilon)$ be a Hopf algebroid.
For  any  left $H$-modules $M_{1}$ and $M_{2}$, $M_{1}\otr M_{2}$ is
again a left $H$-module. Moreover,  the tensor product is
associative:
$(M_{1}\otr M_{2})\otr M_{3}\cong M_{1}\otr (M_{2}\otr M_{3})$.
The category of $H$-modules  in fact becomes  a monoidal category.
\end{pro}

{\bf Example 2.1}  Let $\cald$ denote  the algebra
 of all differential operators
on a smooth manifold $M$, and $R$  the algebra  of smooth
functions on $M$. Then  $\cald$ is a Hopf algebroid over  $R$.
Here, $\alpha =\beta $ is  the inclusion $R\lon \cald $,
while the coproduct $\Delta : \cald \lon \cald \otr \cald$
is defined as
\begin{equation}
\Delta (D)(f, g)=D(fg), \ \  \ \forall D \in \cald , f , \ g\in R.
\end{equation}
Note that  $\cald \otr \cald $  is simply the space of bidifferential
operators. Clearly, $\Delta $ is co-commutative.
 As  for  the co-unit, we take   $\epsilon : \cald \lon R$, 
the natural projection from a differential operator to its $0th$-order part.
In this case, left $\cald$-modules are $D$-modules in the usual sense,
and the tensor product is the usual tensor product
of $D$-modules over   $R$.

We note, however, that this Hopf algebroid does not admit
an antipode in any natural sense.  Given a differential operator
$D$,  its antipode, if
it exists, would be the dual operator $D^*$. However, the latter is
a differential operator on $1$-densities, which
does not possess  any canonical identification with   a differential operator
on $R$.

{\bf Example 2.2} The construction above
 can  be generalized to show that the universal enveloping algebra
$UA$ of a Lie algebroid  $A$ admits a  co-commutative Hopf algebroid
structure. 

Again  we take $R=C^{\infty}(P)$, and let
  $\alpha =\beta : C^{\infty}(M)\lon  UA$ be the inclusion.
For the  co-product, we set
\be
\Delta (f)&=&f\otr 1, \ \  \forall f\in C^{\infty}(P);\\
\Delta (X)&=& X\otr 1+1\otr X, \ \ \forall X\in \gm (A).
\ee
This extends  to a co-product $\Delta : UA\lon
UA\otr UA$ by  the compatibility condition:
Equations (\ref{eq:compatible1}) and (\ref{eq:compatible2}).
Alternatively, we  may identify
$UA$ as the    subalgebra of $\cald (G)$  consisting of  right invariant differential
operators on a (local) Lie groupoid $G$ integrating $A$,
 and then   restrict the co-product $\Delta_{G}$ on $\cald (G)$
to this subalgebra. This is well defined since
$\Delta_G  $ maps  right invariant differential
operators to  right invariant bidifferential
operators. Finally, the co-unit map is defined as
the projection $\epsilon :UA\lon C^{\infty}(P)$.

{\bf Example 2.3} Let $P$ be a smooth manifold and $\cald$ the
ring of differential operators on $P$. Let $\dh$ be
the space of formal power series in $\hbar$ with coefficients
in $\cald$. The Hopf algebroid structure on $\cald$
induces a Hopf algebroid structure on $\dh$, whose structure maps  will
be denoted by the same symbols. 

Let $\phi =1\otr 1+\hbar B_{1}+\cdots \in \cald\otr \cald [[\hbar ]]$
be a  formal power series in $\hbar$ with coefficients being
 bidifferential operators. For any $f, g\in C^{\infty}(P)[[\hbar ]]$, set
$f\hstar g=\phi (f, g)$. 
This product is associative iff the following identity holds:
\begin{equation}
\label{eq:phi}
(\Delta \otr   id)(\phi )\phi^{12}=(id \otr \Delta )(\phi )\phi^{23},
\end{equation}
where $\phi^{12}=\phi \ot 1 \in (\cald \otr \cald )\ot \cald [[\hbar]]$
and $\phi^{23}=1\ot \phi \in \cald \ot (\cald \otr \cald )[[\hbar]]$.
Note that both sides of  the above equation are well defined elements in 
$\cald\otr \cald \otr \cald [[\hbar]]$.  

Equation (\ref{eq:phi})  reminds us  the equation of  a twistor 
defining  a triangular
Hopf algebra (see Section 10 in  \cite{Drinfeld1}).  Thus,
it is not surprising  that our $\phi$ here can be
used to produce  a new Hopf algebroid structure on $\dh$.

Now assume that 
$\phi  =1\otr 1+\hbar B_{1}+\cdots \in \cald\otr \cald [[\hbar ]]$
satisfies Equation (\ref{eq:phi}). Then, $\{f, g\}=B_{1}(f, g)-B_{1}(g, f), \ \ \ \forall
f , \ g\in C^{\infty}(P)$,  is a Poisson bracket,  and
$f\hstar g=\phi (f, g)$ defines a star product on $P$,  which
is a deformation quantization of this  Poisson structure \cite{BEFFL}.

Let $\dhh=\dh$ be  equipped with the usual multiplication,  and
$R_{\hbar}=C^{\infty}(P)[[\hbar ]]$ with the   $*$-product above.
Define $\alpha :R_{\hbar}\lon \dhh$ and $\beta :R_{\hbar}\lon \dhh$,
respectively,  by $\alpha (f)g=f\hstar g$ and $\beta (f)g=g\hstar f$,
$\forall f, g\in C^{\infty}(P)$.
From the associativity of $*_{\hbar}$, it follows that
$\alpha $ is an algebra homomorphism while $\beta$  is an anti-homomorphism,
and $\alpha (R_{\hbar})$ commutes with $\beta (R_{\hbar})$.
Moreover the associativity of $*_{\hbar}$ implies that

\begin{equation}
\label{eq:phi1}
\phi   (\beta (f)\otimes 1-1\ot \alpha (f))=0 \ \mbox{ in } \cald\otr \cald [[\hbar ]], \ \ \forall f\in C^{\infty}(P).
\end{equation}

As an immediate consequence, we have

\begin{lem}
\label{cor:tensor}
Let $M_{1}$ and $M_{2}$ be any $D$-modules. Then the map
\be
\Phi : M_{1} [[\hbar ]] \otrh M_{2}[[\hbar ]]&\lon  &M_{1}\otr M_{2}[[\hbar ]]\\
(m_{1}\oth m_{2})&\lon &\phi \cdot (m_{1}\ot m_{2})
\ee
is  well defined and establishes an isomorphism between
these two vector spaces.
\end{lem}

In particular, when $M_{1}=M_{2}=\cald$, we  obtain
an isomorphism $\Phi : \dhh\otrh \dhh \lon \cald \otr \cald [[\hbar]]$. Now
define  $\Delta_{\hbar}: \dhh \lon \dhh\otrh \dhh $ by
\begin{equation}
\label{eq:coproduct}
\Delta_{\hbar}=\phi^{-1}\Delta \phi ,
\end{equation}
and let  $\epsilon :\dhh\lon R_{\hbar}$  be the projection.
Here, Equation (\ref{eq:coproduct}) means that $\Delta_{\hbar}(x)=\Phi^{-1}(\Delta (x)\phi )$, $\forall x \in \dhh$.

The following result can be easily verified.\\\\
{\bf Theorem A}\ \ {\em $(\dhh , R_{\hbar}, \alpha , \beta , m, \Delta_{\hbar}, \epsilon )$
is a Hopf algebroid.}

Given any $D$-modules $M_{1}$ and $M_{2}$, the 
Hopf algebroid structure on  $\dhh$ induces a $D$-module
structure on $M_{1}[[\hbar ]]\otrh M_{2}[[\hbar ]]$. It is easy
to see that the map $\Phi$ as in Lemma \ref{cor:tensor} is in fact
an isomorphism of $D$-modules.  Hence 
$M_{1}[[\hbar ]]\otrh M_{2}[[\hbar ]]$ and $M_{2}[[\hbar ]]\otrh M_{1}[[\hbar ]]
 $ are isomorphic $D$-modules, where  the isomorphism
is given by $\Phi^{-1}\smalcirc \sigma\smalcirc \Phi$,
with $\sigma$ being  the flipping.


\section{Deformation of Hopf algebroids and quantum groupoids}

\begin{defi}
\label{def:deformation}
A deformation of a Hopf algebroid  $(H, R, \alpha , \beta , m , \Delta , \epsilon)$
over a field $k$ is a topological Hopf algebroid 
$(H_{\hbar}, R_{\hbar}, \alpha_{\hbar} , \beta_{\hbar} , m_{\hbar} , \Delta_{\hbar} , \epsilon_{\hbar})$ over the ring $k[[\hbar ]]$ of 
formal power series in $\hbar $ such that
\begin{enumerate}
\item $H_{\hbar }$ is isomorphic to $H[[\hbar ]]$ as $k[[\hbar ]]$ module with
unit $1_{H}$, and
$R_{\hbar }$ is isomorphic to $R[[\hbar ]]$ as $k[[\hbar ]]$ module with
unit $1_{R}$;
\item  $\alpha_{\hbar}=\alpha (\mbox{mod } \hbar ), 
\beta_{\hbar}=\beta (\mbox{mod } \hbar ), \ \ m_{\hbar}=m  (\mbox{mod } \hbar ),
\epsilon_{\hbar }=\epsilon  (\mbox{mod } \hbar )$;
\item $\Delta_{\hbar}  =\Delta ( \mbox{mod } \hbar ) $
\end{enumerate}
\end{defi}
{\bf Remark} The meaning of (i) and (ii) is clear. However, for 
Condition (iii), we need the following simple fact:
\begin{lem}
Under the  hypothesis (i) and (ii) as in Definition  \ref{def:deformation},
set $V_{\hbar}=H_{\hbar}\otrh H_{\hbar}$. Then $V_{\hbar}/\hbar V_{\hbar}$
is isomorphic to $H\otr H$ as a   vector space.
\end{lem}
 Let  $\tau :H_{\hbar}\otrh H_{\hbar}\lon H\otr H$
denote  the composition of the projection
$V_{\hbar}\lon V_{\hbar}/\hbar V_{\hbar}$ with the isomorphism 
$V_{\hbar}/\hbar V_{\hbar}\cong H\otr H$. We shall also use the
notation $\hbar \mapsto  0$ to denote this map when the
meaning  is clear from the context.
Then, Condition (iii) means that $\limh \Delta_{\hbar} (x)=\Delta (x)$
for any $x\in  \cald$.

\begin{defi}
A quantum groupoid (or  a  QUE algebroid) is a deformation of
the  Hopf algebroid  $UA$ (see Example  2.2) of   a Lie
algebroid $A$.
\end{defi}

Let $(U_{\hbar}A (\cong UA [[\hbar ]]),  R_{\hbar} (\cong C^{\infty}(P)[[\hbar ]]),
 \alpha_{\hbar}, \beta_{\hbar}, m_{\hbar}, \Delta_{\hbar}, \epsilon_{\hbar})$
be a quantum groupoid. It is well known that 
$$\{f , g\}=\limh \frac{1}{\hbar} (f*_{\hbar}g-g*_{\hbar}f), \ \ \
\forall f, g \in C^{\infty}(P)$$
defines a Poisson structure on the base manifold $P$.

Define
\be 
\delta f&=&\limh \frac{1}{\hbar} ( \alpha_{\hbar}f -\beta_{\hbar}f)
\in UA , \ \ \ \forall f\in C^{\infty}(P), \\
\Delta^{1}X&=&\limh \frac{1}{\hbar} (\Delta_{\hbar}X-(1\oth X+X\oth 1))
\in UA \otr UA, 
\ \ \ \forall X\in \gm (A),  \ \mbox{ and }\\
\delta X&=&\Delta^{1}X-\Delta^{1}_{op}X \in UA \otr UA  .
\ee

A routine calculation using the axioms of Hopf algebroids  leads to
\begin{pro}
\label{thm:delta0}
\begin{enumerate}
\item $\delta f\in \gm (A)$ and $\delta X\in \gm (\wedge^{2}A)$ for any
$f \in C^{\infty}(P)$ and $X\in \gm (A); $
\item $\delta (fg)=f\delta g+g\delta f $ for any $f, g \in C^{\infty}(P)$;
\item $\delta (fX)=f\delta X+\delta f\wedge X $ for any $f\in  C^{\infty}(P)$
and $X\in \gm (A)$;
\item $\rho (\delta f)g =\{f , g\}$  for any
$f, g  \in C^{\infty}(P)$, where $\rho :A \lon TP$
is the anchor  of the Lie algebroid $A$.
\end{enumerate}
\end{pro}

Properties (i)-(iii)  allow us to extend $\delta$ to a
well defined degree $1$  derivation  
$\delta :\gm (\wedge^{*}A)\lon \gm (\wedge^{*+1}A)$.
It is not difficult to prove:

\begin{pro}
\label{thm:delta}
\begin{eqnarray}
&&\delta^{2}=0; \ \ \mbox{ and }\\
&&\delta [X, Y]=[\delta X, Y]+[X,  \delta Y ], \ \ \ \ \forall X, Y\in \gm (A).
\label{eq:bialgebroid}
\end{eqnarray}
\end{pro}


\section{Quantization of Lie bialgebroids}

Recall that a Lie bialgebroid  \cite{MackenzieXu1} 
\cite{KS} is a pair of Lie algebroids $(A, A^{*})$
satisfying the compatibility condition: Equation (\ref{eq:bialgebroid}), where
$\delta :\gm (\wedge^{*}A)\lon \gm (\wedge^{*+1}A)$
is the differential induced from the Lie algebroid $A^*$.
Lie bialgebroids include usual Lie bialgebras. Besides, 
for any Poisson manifold $P$,  the pair $(TP , T^{*}P)$
carries a  natural  Lie bialgebroid structure.
More generally, given a  Lie algebroid $A$ and $\Lambda \in \gm (\wedge^{2}A)$
satisfying $[\Lambda ,  \Lambda ]=0$, $(A, A^{*})$   is 
a Lie bialgebroid, called a triangular Lie 
bialgebroid \cite{MackenzieXu1}.
 It is called regular if $\Lambda$ is of constant rank. Another
interesting example arises from the classical dynamical Yang-Baxter equation
(\cite{B}):




{\bf Example 4.1} Let $\frakg$ be 
 a simple Lie algebra  with 
Cartan subalgebra $\frakh$.  A  classical dynamical  r-matrix \cite{EV}
 is a $\frakg \otimes \frakg $
valued  equivariant function $r: \frakh^* \lon \frakg \otimes \frakg$
satisfying  the classical
dynamical Yang-Baxter equation:
$$Alt (dr)+[r^{12}, r^{13}]+[r^{12}, r^{23}]+[r^{13}, r^{23}]=0$$ such that
$r^{12}+r^{21}$ is a constant function  valued
in $(S^{2}\frakg )^{\frakg}$.
Given a  classical dynamical  r-matrix $r$, let $A=T \frakh^{*}\times \frakg$,
equipped with  the product Lie algebroid.
Define $\Lambda \in \gm (\wedge^{2}A)$  by $\Lambda =\sum \xi_{\alpha}\wedge
h_{\alpha}+r$. Here $\{h_{\alpha} \} $ is a basis of $\frakh$,
$\xi_{\alpha}$ its dual basis considered as a  constant 
vector field on $\frakh^{*} $,  and $\xi_{\alpha}\wedge h_{\alpha}$
 is considered as a section of $\wedge^{2}A$
in an evident sense.  It is simple to check that $\Lambda$ 
satisfies the hypothesis as in Theorem 2.1 in  \cite{LX},
and therefore induces a coboundary (or  exact as  called in \cite{LX})
Lie bialgebroid $(A, A^{*})$.

The base space   of  a Lie bialgebroid $(A, A^{*})$ is naturally
equipped with a Poisson structure, which is defined
by the bundle map $\rho \smalcirc \rho_{*}^{*}: T^{*} P\lon TP$.
Here
  $\rho$ and $\rho_{*}$ are  the anchors of $A$ and $A^*$
respectively.

A combination of Propositions \ref{thm:delta0}  and  \ref{thm:delta} leads to \\\\
{\bf Theorem B}\ \  {\em A  quantum groupoid  $U_{\hbar}A$ 
naturally induces a Lie bialgebroid
$(A, A^* )$ as a classical limit. The   induced Poisson
structure  of  this Lie bialgebroid on the base $P$  coincides
with the one induced  from the base  $*$-algebra $R_{\hbar} (\cong 
C^{\infty}(P)[[\hbar ]])$. }

Such a   Lie bialgebroid $(A, A^* )$ is called
the classical limit of the  quantum groupoid  $U_{\hbar}A$.
Conversely,

\begin{defi}
A quantization of  a  Lie bialgebroid  $(A, A^{*})$ is
a quantum groupoid  $U_{\hbar}A$ whose classical limit
is  $(A, A^{*})$.
\end{defi}

It is a deep theorem of Etingof and Kazhdan \cite{EK} that
every Lie bialgebra is quantizable. 
 On the other  hand, the existence of $*$-products for 
arbitrary Poisson manifolds was recently proved  by Kontsevich \cite{K}.
 In terms of Hopf algebroids, this amounts  to saying
that the Lie bialgebroid $(TP, T^{*}P)$  of
a Poisson manifold $P$  is always quantizable.
It is therefore   natural to expect:\\\\
{\bf Conjecture}  Every Lie bialgebroid is quantizable.

In fact, by    modifying Fedosov's method (see \cite{Fedosov} \cite{D}),
 one can prove  the following:\\\\
{\bf Theorem C}\ \ {\em Any regular triangular Lie bialgebroid is quantizable.}


{\bf Remark.}  (1). In  \cite{EV},  Etingof and Varchenko  developed
a theory of $\eta$-Hopf algebroids, which was  aimed  to 
provide   a general  language  for  elliptic quantum groups and
the quantum dynamical Yang-Baxter equation.  The relation
of these studies to  quantization of the
Lie bialgebroids as in   Example 4.1 is the subject of work in progress.

 (2). According to the general spirit of
deformation theory, any deformation  corresponds to
a certain cohomology.
In particular, the deformation of Hopf algebras is controlled
by  the cohomology of a certain  double complex arising from
the Hopf algebra structure \cite{GS}. It is  natural
to ask what is the proper cohomology theory  controlling
the deformation of  a Hopf algebroid, and in particular
what is the premier obstruction to the quantization problem\footnote{The
author is grateful to Drinfeld for raising this question to him.}.

{\bf Acknowledgements.}  The author  would like to thank  Giuseppe Dito,
Vladimir Drinfeld,  Jiang-hua Lu, Masaki Kashiwara, Dale Peterson
 and Alan Weinstein for useful discussions and comments.  Most of
the work was completed during the author's visit in
RIMS, Kyoto. He wishes to thank RIMS, especially  to
his host professor Kyoji  Saito  for his   hospitality
and the support while this project was being done. He is 
also grateful to  Giuseppe Dito for his help in
preparation of  the French abbreviated version.

\end{document}